\documentstyle[11pt]{article}
\begin{document}

\setlength{\leftmargini}{.5\leftmargini}

\newcommand{\ls}[1]
   {\dimen0=\fontdimen6\the\font \lineskip=#1\dimen0
\advance\lineskip.5\fontdimen5\the\font \advance\lineskip-\dimen0
\lineskiplimit=.9\lineskip \baselineskip=\lineskip
\advance\baselineskip\dimen0 \normallineskip\lineskip
\normallineskiplimit\lineskiplimit \normalbaselineskip\baselineskip
\ignorespaces }


\def\thepart{\Roman{part}} 
\def\thesection {\arabic{section}.}
\def\thesubsection {\thesection\arabic{subsection}.}
\def\thesubsubsection {\thesubsection\arabic{subsubsection}.}
\def\theparagraph {\thesubsubsection\arabic{paragraph}.}
\def\thesubparagraph {\theparagraph\arabic{subparagraph}.}

    \setcounter{equation}{0}

\newcommand{\hsp}{{\hspace*{\parindent}}}

\newtheorem{problem}{Problem}[section]
\newtheorem{definition}{Definition}[section]
\newtheorem{lemma}{Lemma}
\newtheorem{proposition}{Proposition}
\newtheorem{corollary}{Corollary}[section]
\newtheorem{example}{Example}[section]
\newtheorem{conjecture}{Conjecture}
\newtheorem{algorithm}{Algorithm}[section]
\newtheorem{theorem}{Theorem}
\newtheorem{exercise}{Exercise}[section]

\makeatletter
\def\@begintheorem#1#2{\it \trivlist \item[\hskip \labelsep{\bf #1\
#2.}]}
\makeatother
                                                   
\def\P{{\bf P}}
\def\E{{\bf E}}
\newcommand{\be}{\begin{equation}}
\newcommand{\ee}{\end{equation}}
\newcommand{\bea}{\begin{eqnarray}}
\newcommand{\eea}{\end{eqnarray}}

\newcommand{\beq}[1]{\begin{equation}\label{#1}}
\newcommand{\eeq}{\end{equation}}
\newcommand{\req}[1]{(\ref{#1})}
\newcommand{\beqn}[1]{\begin{eqnarray}\label{#1}}
\newcommand{\eeqn}{\end{eqnarray}}

\newcommand{\beaa}{\begin{eqnarray*}}
\newcommand{\eeaa}{\end{eqnarray*}}

\newcommand{\eq}[1]{(\ref{#1})}

\def\le{\leq}
\def\ge{\geq}
\def\lt{<}
\def\gt{>}

\newcommand{\lip}{\langle}
\newcommand{\rip}{\rangle}
\newcommand{\uu}{\underline}
\newcommand{\oo}{\overline}
\newcommand{\La}{\Lambda}
\newcommand{\la}{\lambda}
\newcommand{\eps}{\varepsilon}
\newcommand{\vp}{\varphi}

\newcommand{\dint}{\displaystyle\int}
\newcommand{\dsum}{\displaystyle\sum}
\newcommand{\dfr}{\displaystyle\frac}
\newcommand{\bige}{\mbox{\Large\it e}}
\newcommand{\integers}{Z\!\!\!Z}
\newcommand{\rationals}{{\rm I\!Q}}
\newcommand{\reals}{{\rm I\!R}}
\newcommand{\realsd}{\reals^d}
\newcommand{\NN}{{\rm I\!\!N}}
\newcommand{\degree}{{\scriptscriptstyle \circ }}
\newcommand{\dfn}{\stackrel{\triangle}{=}}
\def\complex{\mathop{\raise .45ex\hbox{${\bf\scriptstyle{|}}$}
     \kern -0.40em {\rm \textstyle{C}}}\nolimits}
\def\hilbert{\mathop{\raise .21ex\hbox{$\bigcirc$}}\kern -1.005em {\rm\textstyle{H}}} 

\newcommand{\calA}{{\cal A}}
\newcommand{\calC}{{\cal C}}
\newcommand{\calD}{{\cal D}}
\newcommand{\calF}{{\cal F}}
\newcommand{\calL}{{\cal L}}
\newcommand{\calM}{{\cal M}}
\newcommand{\calP}{{\cal P}}
\newcommand{\calX}{{\cal X}}

\newcommand{\Prob}{{\rm Prob\,}}
\newcommand{\mod}{{\rm mod\,}}
\newcommand{\sinc}{{\rm sinc\,}}
\newcommand{\ctg}{{\rm ctg\,}}
\newcommand{\ifff}{\mbox{\ if and only if\ }}
\newcommand{\proof}{\noindent {\bf Proof:\ }}
\newcommand{\remark}{\noindent {\bf Remark:\ }}
\newcommand{\remarks}{\noindent {\bf Remarks:\ }}
\newcommand{\note}{\noindent {\bf Note:\ }}

\newcommand{\boldx}{{\bf x}}
\newcommand{\boldX}{{\bf X}}
\newcommand{\boldy}{{\bf y}}
\newcommand{\uux}{\uu{x}}
\newcommand{\uuY}{\uu{Y}}

\newcommand{\liml}{\underline{\lim}_{l\rightarrow \infty}}
\newcommand{\limn}{\underline{\lim}_{n \rightarrow \infty}}
\newcommand{\limN}{\lim_{N \rightarrow \infty}}
\newcommand{\limr}{\lim_{r \rightarrow \infty}}
\newcommand{\limd}{\lim_{\delta \rightarrow \infty}}
\newcommand{\limM}{\lim_{M \rightarrow \infty}}
\newcommand{\limsupn}{\limsup_{n \rightarrow \infty}}
\newcommand{\liminfn}{\liminf_{n \rightarrow \infty}}

\newcommand{\imii}{\int_{-\infty}^{\infty}}
\newcommand{\imix}{\int_{-\infty}^x}
\newcommand{\ioi}{\int_o^\infty}

\newcommand{\ARROW}[1]
  {\begin{array}[t]{c}  \longrightarrow \\[-0.2cm] \textstyle{#1} \end{array} }

\newcommand{\AR}
 {\begin{array}[t]{c}
  \longrightarrow \\[-0.3cm]
  \scriptstyle {j\rightarrow \infty}
  \end{array}}

\newcommand{\ARn}
 {\begin{array}[t]{c}
  \longrightarrow \\[-0.3cm]
  \scriptstyle {n\rightarrow \infty}
  \end{array}}
\newcommand{\pile}[2]
  {\left( \begin{array}{c}  {#1}\\[-0.2cm] {#2} \end{array} \right) }

\newcommand{\floor}[1]{\left\lfloor #1 \right\rfloor}

\newcommand{\mmbox}[1]{\mbox{\scriptsize{#1}}}

\newcommand{\ffrac}[2]
  {\left( \frac{#1}{#2} \right)}

\newcommand{\one}{\frac{1}{n}\:}
\newcommand{\half}{\frac{1}{2}\:}
\newcommand{\cE}{{\cal E}}
\newcommand{\cP}{{\cal P}}
\newcommand{\bY}{{\bf Y}}
\newcommand{\bX}{{\bf X}}
\newcommand{\bZ}{{\bf Z}}
\newcommand{\s}{{\sigma}}

\def\squarebox#1{\hbox to #1{\hfill\vbox to #1{\vfill}}}
\newcommand{\qed}{\hspace*{\fill}
           \vbox{\hrule\hbox{\vrule\squarebox{.667em}\vrule}\hrule}\smallskip}

\ls{1.5}

\begin{centering}
{\large \bf A conjecture concerning
optimality of the Karhunen-Loeve basis in nonlinear reconstruction 
}
\\[2em]
{\sc Stephane Mallat},\footnote{Ecole Polytechnique, Palaiseau}
and
{\sc Ofer Zeitouni}
\footnote{Technion, Haifa; currently, Weizmann Institute and University 
of Minnesota}

September 4, 2000. This version September 15, 2011.

\end{centering}

\section{Introduction}
The problem posed in this note has its root in discussion carried out
more than 10 years ago between the authors. Since then, we have discussed it 
with numerous people, and it has been posted as an open problem on the web site
of one of us. We decided to post it on the ArXiv in order to have a permanent
and stable version for it.
 
\noindent
\section{The  conjecture}
Let $\bY:=(Y_1,\ldots,Y_N)$ denote an $N$-dimensional
 Gaussian vector with independent zero mean components of variance 
$\sigma_i=E(Y_i^2)$. We assume for concreteness that $\sigma_i\geq\sigma_{i+1}$.

Let $T(\theta)$ denote an arbitrary orthogonal matrix on $\reals^N$ ($\theta$
is a $N(N-1)/2$ dimensional parameter, and we take $T(0)=I$), 
and define the random variable
$\bX(\theta)=T(\theta) \bY$. For any $M<N$, define
$$ \cE(N,M,\theta)=E(\min_{i_1\neq i_2\neq \ldots\neq  i_{N-M}} \sum_{j=1}^{N-M}
(X_{i_j}(\theta))^2)\,.$$
$\cE(N,M,\theta)$ is the mean square error when reconstructing $\bX(\theta)$ 
according to 
its $M$ largest (in absolute value) components. When the matrix $T(\theta)$ is
a permutation matrix,  this is the reconstruction error when keeping the
$M$ largest (in absolute value) components of $\bY$, where $\bY$ is already
expressed in its Karhunen-Loeve basis, whereas other choices of $\theta$
correspond to an expansion in other, non K-L bases.

Let $T(\theta)\in \cP$ if $T(\theta)$ is composed only of zeroes and ones, i.e.
$T(\theta)$ is a permutation and reflection matrix.
We have the following
\begin{conjecture}
\label{conj-1}
$$\min_{\theta}\cE(N,M,\theta)= \min_{T(\theta)\in \cP} \cE(N,M,\theta)=
\cE(N,M,0)\,.$$
\end{conjecture}

Conjecture \ref{conj-1}, if true, implies that the Karhunen-Loeve 
basis is the best basis not only for linear reconstruction but also 
for nonlinear reconstruction based on the $M$ largest projections.
 
\section{The case of $M=1$}
We do not know how to prove in general 
Conjecture \ref{conj-1}. However, it does hold true for $M=1$, i.e
reconstruction based on the largest projection. Indeed, we have
\begin{lemma}
\label{lem-1}
Conjecture \ref{conj-1} holds true if  $M=1$.
\end{lemma}

\proof
We can re-parametrize $T$ such that $\eta_i=E(X_i^2)$ satisfies
$\eta_i\geq \eta_{i+1}$. Note that $\eta_i=\sum_{j=1}^N T_{ij}^2 \sigma_j$.
Clearly, it is enough then to prove that 
$E(\max X_i^2)\leq E(\max Y_i^2)$.

Let $\tilde  \bX$ denote a vector of independent Gaussian random variables
with $E\tilde X_i=0$ and $E(\tilde X_i)^2=\eta_i$. By Sidak's inequality
\cite{6aut},
for any $t>0$,
$$P(\cap \{|X_i|<t\})\geq \prod P(|X_i|<t) = P(\cap \{|\tilde X_i|<t\})\,,$$
implying that
\be
\label{cor-ineq}
 E(\max |X_i|^2)\leq E(\max |\tilde X_i|^2)\,.
\ee

Next, we can check that 
\be
\label{shur}
\sum_{i=1}^m \sigma_i \geq \sum_{i=1}^m \eta_i,\; i=1,\ldots , N-1\,.
\ee
and clearly, because $\sum_i T_{ij}^2=1$, also 
$$\sum_{i=1}^N \sigma_i=\sum_{i=1}^N \eta_i\,.$$ 
Indeed, when $m=1$, (\ref{shur}) holds because $\sum_j T_{1j}^2 =1$.
For $m=2$, we have that
\begin{eqnarray*}
\eta_1+\eta_2&=& \sum_j (T_{1j}^2+T_{2j}^2) \sigma_j\\
&=& 
\sum_{j=2}^N 
 (T_{1j}^2+T_{2j}^2) \sigma_j + (T_{11}^2+T_{21}^2 )\sigma_1\\
&\leq &
\sigma_2 \sum_{j=2}^N 
 (T_{1j}^2+T_{2j}^2)  
+ (T_{11}^2+T_{21}^2 )\sigma_1\\
&=& \sigma_2(2-T_{11}^2-T_{21}^2)
+ (T_{11}^2+T_{21}^2 )\sigma_1\\
&=& \sigma_1+\sigma_2 +(\sigma_2-\sigma_1)(1-T_{11}^2-T_{21}^2)
\leq \sigma_1+\sigma_2\,.
\end{eqnarray*}
The general case of (\ref{shur}) follows by induction.

By an inequality of Marshall and Proschan, see \cite[Application 7.A.18]{SS},
one concludes that
for any convex, permutation symmetric function $\phi$,
$$ E(\phi(\tilde X_1^2,\ldots,\tilde X_N^2)) \leq
 E(\phi( Y_1^2,\ldots, Y_N^2))\,.$$
Applying this to the function $\phi(\cdot)=\max x_i^2$, one concludes that
$E(\max Y_i^2)\geq E(\max (\tilde X_i)^2)$, which together with
(\ref{cor-ineq}) yields 
$E(\max Y_i^2)\geq E(\max X_i^2)$, as claimed.
\qed

Remark: the Schur convexity part of the argument holds also for the function
$$\phi_M({\bf x})=
\max_{i_1\neq i_2\neq \ldots\neq i_M} \sum_{j=1}^M x_{i_j}^2\,.$$
What is missing in order to prove the conjecture for general $M$
is the analog of (\ref{cor-ineq}): is it true that
\begin{equation}
\label{eq-11}
 E(\phi_M({\bf X}))\leq
 E(\phi_M(\tilde {\bf X}))\,?
\end{equation}
{\it Added September 15, 2011:
R. van Handel communicated to us
the following  counter example to (\ref{eq-11}): take $N=3$, $M=2$
and $Z_1,Z_2,Z_3$ three independent standard Gaussians. Define
$X_1=(Z_1-Z_2)/\sqrt{2}$, $X_2=(Z_2-Z_3)/\sqrt{2}$, and
$X_3=(Z_3-Z_1)/\sqrt{2}$. The corresponding $\tilde X_i$ are independent
standard Gaussians. One checks numerically that
$0.17\sim E\min(X_i^2)<E \min(\tilde X_i^2)\sim 0.19$ 
(note that Conjecture \ref{conj-1}
does hold in this case). Of course, it is possible that
using (\ref{eq-11}) for only a subset of all $T$s can help.

This example also disproves the conjecture in \cite[Problem 6, pg. 279]{lifshitz}}

\noindent
{\bf Remark} Some inequalities related to the problem discussed in this
note can be found in \cite{GL}. However, the results contained there
are not enough to resolve Conjecture \ref{conj-1}, even within 
a multiplicative factor.

\bibliographystyle{plain}

\end{document}